\documentclass[11pt]{amsart}
\usepackage{amssymb, latexsym}
\theoremstyle{plain}
\newtheorem{theorem}{Theorem}

\newtheorem*{thm-cheb}{Theorem (Chebyshev)}

\newtheorem*{PI}{Theorem P}

\newtheorem*{2'}{Theorem 2'}
\newtheorem*{3'}{Theorem 3'}

\theoremstyle{remark}

\newtheorem*{Remark 1}{Remark 1}
\newtheorem*{Remark 2}{Remark 2}
\newtheorem*{Remark 3}{Remark 3}
\newtheorem*{Remark 4}{Remark 4}

\numberwithin{equation}{section}

\begin{document}

\title[Secretary problem with biased arrival  via  Mallows distributions]
 {The secretary problem with biased arrival order via a Mallows distribution}

\author{Ross G. Pinsky}


\address{Department of Mathematics\\
Technion---Israel Institute of Technology\\
Haifa, 32000\\ Israel}
\email{ pinsky@math.technion.ac.il}

\urladdr{https://pinsky.net.technion.ac.il/}

\subjclass[2000]{60G40, 60C05} \keywords{secretary problem, optimal stopping, Mallows distribution, inversions}
\date{}

\begin{abstract}
We solve the secretary problem in the case that the ranked items arrive in a statistically biased order rather than in uniformly random order.
The bias is given by a Mallows distribution with parameter $q\in(0,1)$, so that higher ranked items tend to arrive later and lower ranked items
tend to arrive sooner.
In the classical problem, the asymptotically optimal strategy is to reject the first  $M_n^*$  items, where $M_n^*\sim\frac ne$, and then to select
the first item ranked higher than any of the first  $M_n^*$ items (if such an item exists). This yields $\frac1e$ as the  limiting   probability of success.
The Mallows distribution with parameter $q=1$ is the uniform distribution.
For the regime $q_n=1-\frac cn$, with $c>0$, the case of weak bias, the optimal strategy occurs with
$M_n^*\sim n\Big(\frac1c\log\big(1+\frac{e^c-1}e\big)\Big)$, with the limiting probability of success being $\frac1e$.
For the regime $q_n=1-\frac c{n^\alpha}$, with $c>0$ and $\alpha\in(0,1)$, the case of moderate bias,
the optimal strategy occurs with
$n-M_n\sim\frac{n^\alpha}c$, with the limiting probability of success being $\frac1e$. For  fixed $q\in(0,1)$, the case of strong bias,
the optimal strategy occurs with
$M_n^*=n-L$ where $\frac{L-1}L<q\le \frac L{L+1}$, with limiting probability of success being $(1-q)q^{L-1}L>\frac1e$.
\end{abstract}

\maketitle

\section{Introduction and Statement of  Results}
Recall the classical secretary problem: For $n\in\mathbb{N}$, a set of $n$ ranked items is revealed, one item at a time, to an observer whose objective is to select the item with
the highest rank. The order of the items is completely random; that is, each of the  $n!$ permutations of the ranks is equally likely.
At each stage, the observer only knows the relative ranks of the items that have arrived thus far, and must  either select the current item, in which case the process terminates, or reject it and continue to the next item. If the observer rejects the first $n-1$ items, then the $n$th and final item to arrive must be accepted.
As is very well known, asymptotically as $n\to\infty$, the optimal  strategy is to reject the first $M_n$ items, where $M_n\sim \frac ne$, and then to select the first
later-arriving item whose rank is higher than that of any of the first $M_n$ items (if such an item exists).
The limiting probability of successfully selecting the item of highest rank is $\frac1e$.

Over the years, the secretary problem has been generalized in many directions. See the 1989 paper \cite{F} for a history of the problem and some natural generalizations.
Many of the more recent papers concerning the secretary problem are in the computer science literature.

In this paper, we consider the secretary problem in the case that the order of arrival is biased
so that there is a statistical tendency for higher ranked items to arrive later and lower ranked items to arrive sooner (or vice versa).
We are only aware of one paper in the literature that considers such a situation. The paper gives sufficient conditions on the permutation distribution to guarantee
that there exists an algorithm for which the probability of success is bounded away from zero, independent of $n$
\cite{KKN}.

 The bias we introduce is  via a  Mallows distribution on the set $S_n$ of permutations of $[n]$. \it Our convention will be that the number $n$ represents the highest ranking and the number 1 represents the lowest
ranking.\rm\
Recall that the Mallows distribution with parameter $q>0$ is the
distribution $P_n^q$ for which  $P_n^q(\sigma)$ is proportional to $q^{\text{inv}(\sigma)}$, for $\sigma\in S_n$, where
$\text{inv}(\sigma)$ is equal to the number of inversions in the permutation $\sigma$.
Thus, for $q\in(0,1)$, $P_n^q(\sigma)$ is decreasing in $\text{inv}(\sigma)$, while
for $q>1$, it is increasing in $\text{inv}(\sigma)$. Therefore, larger numbers have a tendency to appear toward the end of the permutation
if $q\in(0,1)$, and toward the beginning of the permutation if $q>1$. Of course, $q=1$ corresponds to the uniform distribution on $S_n$.
Recall that the reverse  of a permutation $\sigma=\sigma_1\cdots\sigma_n$ is the permutation $\sigma^{\text{rev}}:=\sigma_n\cdots\sigma_1$.
The Mallows distributions satisfy the following duality between $q>1$ and $q<1$:
$$
P_n^q(\sigma)=P_n^{\frac1q}(\sigma^{\text{rev}}),\ \text{for}\ q>0, \sigma\in S_n\ \text{and}\ n=1,2,\cdots.
$$
Consequently, it suffices to restrict our study to $q\in (0,1)$.
In terms of the secretary problem, this means that there is a tendency for the higher ranked items to arrive later and the lower ranked items to arrive sooner.

Define $\mathcal{I}_n(\sigma)=\text{inv}(\sigma)$, for $\sigma\in S_n$. It is well known that
under the uniform distribution on $S_n$, the random variable $\mathcal{I}_n$ satisfies the weak law of large numbers in the form
$\text{\rm w-}\lim_{n\to\infty}\frac{\mathcal{I}_n}{n^2}=\lim_{n\to\infty}\frac{E_n^1\mathcal{I}_n}{n^2}=\frac14$.
In \cite{Pin} the behavior of $\mathcal{I}_n$  under Mallows distributions was investigated and the following results were proven.
\begin{PI}\label{PI}
\noindent
\noindent i. Under $P_n^{q_n}$,  for $q_n=1-\frac cn$, with $c>0$,
$$
\text{\rm w-}\lim_{n\to\infty}\frac{\mathcal{I}_n}{n^2}=
\lim_{n\to\infty}\frac{E_n^{q_n}\mathcal{I}_n}{n^2}=
\frac1{c^2}\int_0^{1-e^{-c}}\big(\frac1{1-x}+\frac{\log(1-x)}x\big)dx:=I(c);
$$
furthermore, $\lim_{c\to\infty}I(c)=0$ and $\lim_{c\to0}I(c)=\frac14$;

\noindent ii. Under $P_n^{q_n}$,  for $q_n=1-\frac c{n^\alpha}$, with $c>0$ and $\alpha\in(0,1)$,
$$
\text{\rm w-}\lim_{n\to\infty}\frac{\mathcal{I}_n}{n^{1+\alpha}}=
\lim_{n\to\infty}\frac{E_n^{q_n}\mathcal{I}_n}{n^{1+\alpha}}=
\frac1c;
$$

iii. Under $P_n^q$,  for $q\in(0,1)$,
$$
\text{\rm w-}\lim_{n\to\infty}\frac{\mathcal{I}_n}n=\lim_{n\to\infty}\frac{E_n^q\mathcal{I}_n}n=\frac q{1-q}.
$$
\end{PI}

In light of the above result, we will say that  the sequence of distributions $\{P_n^{q_n}\}_{n=1}^\infty$
corresponds to weak, moderate or strong bias respectively according to whether  $\{q_n\}_{n=1}^\infty$ is as in part (i), (ii) or (iii) of  Theorem P.
Returning to the secretary problem, for each $n\in \mathbb{N}$, let $\mathcal{S}(n,M)$, $0\le M\le n-1$, denote the strategy whereby the observer rejects the first $M$ items, and then
 selects the first
later-arriving item whose rank is higher than that of any of the first $M$ items (if such an item exists).
If the order of arrival of the items is biased via the Mallows distribution with parameter $q\in(0,1)$,
let $\mathcal{P}_n^q(\mathcal{S}(n,M))$ denote the probability of successfully selecting the item of highest rank.
The following theorem determines the asymptotically optimal strategy $\mathcal{S}(n,M_n^*)$  and the limiting optimal probability of success for the weak, the moderate and the strong regimes of bias.
In particular, the limiting optimal probability is $\frac1e$ in the cases of weak or moderate bias, but is larger in the case of strong bias.

\begin{theorem}\label{SecMall}

\noindent i. Let $q_n=1-\frac cn$, where $c>0$. Then the asymptotically optimal strategy is $\mathcal{S}(n,M_n^*)$, where
\begin{equation}\label{optimali}
M_n^*\sim n\Big(\frac1c\log\big(1+\frac{e^c-1}e\big)\Big),
\end{equation}
 and the corresponding  limiting probability of success
is $\frac1e$:
$$
\lim_{n\to\infty}\mathcal{P}_n^{q_n}(\mathcal{S}(n,M_n^*))=\frac1e.
$$
Also, $\lim_{c\to\infty} \frac1c\log\big(1+\frac{e^c-1}e\big)=1$ and $\lim_{c\to0} \frac1c\log\big(1+\frac{e^c-1}e\big)=\frac1e$.

\noindent ii.
Let $q_n=1-\frac c{n^\alpha}$, where $c>0$ and $\alpha\in(0,1)$. Then the asymptotically optimal strategy  is $\mathcal{S}(n,M_n^*)$, where
$$
n-M_n^*\sim\frac{n^\alpha}c,
$$
 and the corresponding limiting probability of success
is $\frac1e$:
$$
\lim_{n\to\infty}\mathcal{P}_n^{q_n}(\mathcal{S}(n,M_n^*))=\frac1e.
$$
\noindent iii. Let $q\in(0,1)$.  Then the asymptotically optimal strategy  is $\mathcal{S}(n,M_n^*)$, where
$$
M_n^*=n-L, \ \text{if}\ \  \frac {L-1}L<q\le \frac L{L+1},\ L=1,2,\cdots.
$$
The corresponding limiting probability of success is given by
$$
\lim_{n\to\infty}\mathcal{P}_n^q(\mathcal{S}(n,M_n^*))=(1-q)q^{L-1}L,\ \text{if}\   \frac {L-1}L<q\le \frac L{L+1},\ L=1,2,\cdots.
$$
In particular,
$$
\lim_{n\to\infty}\mathcal{P}_n^q(\mathcal{S}(n,M_n^*))>\frac1e, \ 0<q<1;
$$
$$
\lim_{q\to1}\lim_{n\to\infty}\mathcal{P}_n^q(\mathcal{S}(n,M_n^*))=\frac1e,
$$
and
$$
\lim_{n\to\infty}\mathcal{P}_n^q(\mathcal{S}(n,M_n^*))=1-q, \ 0< q\le \frac12.
$$
\end{theorem}

The following theorem gives the exact behavior of $\mathcal{P}_n^q(\mathcal{S}(n,M))$ for any $n,q,M$.
\begin{theorem}\label{Exact}
For  $n\in\mathbb{N}, M\in\{0,1,\cdots, n-1\}$ and $q\in(0,1)$,
\begin{equation}\label{exactform}
\mathcal{P}_n^q(\mathcal{S}(n,M))=\begin{cases}\frac{1-q}{1-q^n}q^{n-M-1}(1-q^M)\sum_{j=M+1}^n\frac1{1-q^{j-1}},\ \text{if}\ M\ge1;\\
\frac{1-q}{1-q^n}q^{n-1},\ \text{if}\ M=0.\end{cases}
\end{equation}
\end{theorem}

We prove Theorem \ref{Exact} in section \ref{Exactpf} and then use this result to prove  Theorem \ref{SecMall} in section \ref{SecMallpf}.

\section{Proof of Theorem \ref{Exact}}\label{Exactpf}

 For the proof of the theorem, will need the following so-called online construction of a random permutation in $S_n$ distributed according to the Mallows distribution with parameter $q$. By ``online'' we mean that
the random permutation is constructed in $n$ steps, with one number being added to the permutation at each step.
Let $\{X_j\}_{j=2}^n$ be independent random variables with $X_j$ distributed as a geometric random variable with parameter $1-q$ and truncated at $j-1$; that is,
\begin{equation}\label{Xdist}
P(X_j=m)=\frac{(1-q)q^m}{1-q^j},\ m=0,\cdots, j-1.
\end{equation}
Consider a horizontal line on which to place the numbers in $[n]$.
We begin by placing down the number 1. Then inductively, if  we have already placed down the numbers $1,2,\cdots, j-1$, the number $j$ gets placed down in the position for which there are $X_j$ numbers to its right.
Thus, for example, for $n=4$, if $X_2=1$, $X_3=2$ and $X_4=0$, then we obtain the permutation 3214.
To see that this construction does indeed induce the Mallows distribution with parameter $q$,
 note that the number of inversions in the constructed permutation $\sigma$ is $\sum_{j=2}^nX_j$. Thus, letting
 $Z_n(q)=\prod_{k=2}^n\frac{1-q^k}{1-q}$   and  using \eqref{Xdist}, one obtains
 $P(X_j=x_j,\ j=2,\cdots, n)=\frac1{Z_n(q)}q^{\sum_{j=2}^nx_j}=\frac{q^{\text{inv}(\sigma)}}{Z_n(q)}$.
This calculation also yields the normalization constant $Z_n(q)$.

The following fact will also be essential in the proof of the theorem. One has $\text{inv}(\sigma)=\text{inv}(\sigma^{-1})$, for all $\sigma\in S_n$, where
$\sigma^{-1}$ denotes the inverse permutation of the permutation $\sigma$. Therefore, if $\sigma\in S_n$ is distributed according to the Mallows distribution
with parameter $q$, then $\sigma^{-1}$ also has this distribution.

Let $\sigma=\sigma_1\sigma_2\cdots\sigma_n\in S_n$ represent the rankings of the $n$ items that arrive one by one. That is, $\sigma_j$ is the ranking of the $j$th item to arrive. First consider the case $M=0$. The strategy $\mathcal{S}(n,0)$ will select the highest ranking item if and only if $\sigma_1=n$. From the online construction above,
$P_n^q(\sigma_1=n)=P(X_n=n-1)=\frac{1-q}{1-q^n}q^{n-1}$. This gives \eqref{exactform} for the case $M=0$.

From now on, assume that $M\ge1$.
Then the strategy
$\mathcal{S}(n,M)$ will select the highest ranking item if and only if for some $j\in\{M+1,\cdots, n\}$, one has
$\sigma_j=n$ and $\max(\sigma_1,\cdots,\sigma_{j-1})=\max(\sigma_1,\cdots, \sigma_M)$.
So
\begin{equation}\label{exactstep1}
\mathcal{P}_n^q(\mathcal{S}(n,M))=\sum_{j=M+1}^nP_n^q( \sigma_j=n,\ \max(\sigma_1,\cdots,\sigma_{j-1})=\max(\sigma_1,\cdots, \sigma_M)).
\end{equation}
Since $\sigma$ and $\sigma^{-1}$ have the same distribution under $P_n^q$, we have
\begin{equation}\label{useinverse}
\begin{aligned}
&P_n^q(\sigma_j=n,\ \max(\sigma_1,\cdots,\sigma_{j-1})=\max(\sigma_1,\cdots, \sigma_M))=\\
&P_n^q(\sigma_j^{-1}=n,\ \max(\sigma^{-1}_1,\cdots,\sigma^{-1}_{j-1})=\max(\sigma^{-1}_1,\cdots, \sigma^{-1}_M)).
\end{aligned}
\end{equation}
From the online construction, the events
$\{\sigma_j^{-1}=n\}$ and $\{\max(\sigma^{-1}_1,\cdots,\sigma^{-1}_{j-1})=\max(\sigma^{-1}_1,\cdots, \sigma^{-1}_M)\}$ are independent.
Indeed, the former event depends only on $\{X_l\}_{l=j}^n$ and the latter event depends only on
$\{X_l\}_{l=M+1}^{j-1}$.
Thus, from \eqref{useinverse}, 
\begin{equation}\label{separate}
\begin{aligned}
&P_n^q(\sigma_j=n,\ \max(\sigma_1,\cdots,\sigma_{j-1})=\max(\sigma_1,\cdots, \sigma_M))=\\
&P_n^q(\sigma_j^{-1}=n)P_n^q(\max(\sigma^{-1}_1,\cdots,\sigma^{-1}_{j-1})=\max(\sigma^{-1}_1,\cdots, \sigma^{-1}_M)).
\end{aligned}
\end{equation}
Again using the fact that $\sigma$ and $\sigma^{-1}$ have the same distribution under $P_n^q$, 
we have
\begin{equation}\label{easy}
P_n^q(\sigma^{-1}_j=n)=P_n^q(\sigma_j=n)=P(X_n=n-j)=\frac{(1-q)q^{n-j}}{1-q^n}.
\end{equation}
Also, from the online construction,
\begin{equation}\label{maxpart}
\begin{aligned}
&P_n^q(\max(\sigma^{-1}_1,\cdots,\sigma^{-1}_{j-1})=\max(\sigma^{-1}_1,\cdots, \sigma^{-1}_M))=P(X_l\ge1,\ l=M+1,\cdots, j-1)=\\
&\prod_{l=M+1}^{j-1}\big(1-\frac{1-q}{1-q^l}\big)=
\frac{q^{j-1-M}\big(\prod_{l=M+1}^{j-1}(1-q^{l-1})\big)}{\prod_{l=M+1}^{j-1}(1-q^l)}=q^{j-1-M}\frac{1-q^M}{1-q^{j-1}}.
\end{aligned}
\end{equation}
Now \eqref{exactstep1} and \eqref{separate}-\eqref{maxpart} yield
\begin{equation}
\begin{aligned}
&\mathcal{P}_n^q(\mathcal{S}(n,M))=\sum_{j=M+1}^n\frac{(1-q)q^{n-j}}{1-q^n}q^{j-1-M}\frac{1-q^M}{1-q^{j-1}}=\\
&\frac{1-q}{1-q^n}q^{n-M-1}(1-q^M)\sum_{j=M+1}^n\frac1{1-q^{j-1}}.
\end{aligned}
\end{equation}
\hfill$\square$

\section{Proof of Theorem \ref{SecMall}}\label{SecMallpf}
\noindent \it Proof of part (i).\rm\
Let $c>0$ be fixed and let  $b\in(0,1)$ vary. Substituting $q=q_n=1-\frac cn$ and $M=M_n\sim bn$ in \eqref{exactform}, we obtain
\begin{equation}\label{beginasymp}
\begin{aligned}
&\mathcal{P}_n^{q_n}(\mathcal{S}(n,M_n))\sim \frac cn\frac1{1-(1-\frac cn)^n}(1-\frac cn)^{(1-b)n}(1-(1-\frac cn)^{bn})\sum_{j=[bn]+1}^n\frac1{1-(1-\frac cn)^{j-1}}\sim\\
&\frac {ce^{-c(1-b)}(1-e^{-cb})}{1-e^{-c}}\thinspace\frac1n\sum_{j=[bn]+1}^n\frac1{1-(1-\frac cn)^{j-1}}=\frac{ce^{-c}(e^{bc}-1)}{1-e^{-c}}\frac1n\sum_{j=[bn]+1}^n\frac1{1-(1-\frac cn)^{j-1}}.
\end{aligned}
\end{equation}
We have
$$
\frac1n\sum_{j=[bn]+1}^n\frac1{1-(1-\frac cn)^{j-1}}\sim\frac1n\sum_{j=[bn]+1}^n\frac1{1-e^{-\frac{c(j-1)}n}},
$$
and thus,
\begin{equation}\label{intlim}
\lim_{n\to\infty}\frac1n\sum_{j=[bn]+1}^n\frac1{1-(1-\frac cn)^{j-1}}=\int_b^1\frac1{1-e^{-cx}}dx.
\end{equation}
Making the substitution, $y=e^{-cx}$, we have
\begin{equation}\label{intcalc}
\begin{aligned}
&\int_b^1\frac1{1-e^{-cx}}dx=\frac1c\int_{e^{-c}}^{e^{-bc}}\frac1{y(1-y)}dy=
\frac1c\int_{e^{-c}}^{e^{-bc}}(\frac1y+\frac1{1-y})dy=\\
&1-b+\frac1c\log\frac{1-e^{-c}}{1-e^{-bc}}.
\end{aligned}
\end{equation}
From \eqref{beginasymp}-\eqref{intcalc}, we obtain
\begin{equation}\label{limprobi}
\lim_{n\to\infty}\mathcal{P}_n^{q_n}(\mathcal{S}(n,M_n))=\frac{ce^{-c}(e^{bc}-1)}{1-e^{-c}}\big(1-b+\frac1c\log\frac{1-e^{-c}}{1-e^{-bc}}\big),
\ \text{if}\ M_n\sim bn.
\end{equation}

Define
\begin{equation}\label{H}
H(b)=(e^{bc}-1)\big(1-b+\frac1c\log\frac{1-e^{-c}}{1-e^{-bc}}\big).
\end{equation}
Note that $H(0^+)=H(1^-)=0$.
We have
\begin{equation}
\begin{aligned}
&H'(b)=e^{bc}c(1-b+\frac1c\log\frac{1-e^{-c}}{1-e^{-bc}})+(e^{bc}-1)(-1-\frac{e^{-bc}}{1-e^{-bc}})=\\
&e^{bc}c(1-b+\frac1c\log\frac{1-e^{-c}}{1-e^{-bc}})-e^{bc}.
\end{aligned}
\end{equation}
Thus, $H'(b)=0$ if and only if
$$
c(1-b+\frac1c\log\frac{1-e^{-c}}{1-e^{-bc}})=1,
$$
or equivalently,
$$
e^{1-(1-b)c}=\frac{1-e^{-c}}{1-e^{-bc}}.
$$
Solving for $b$ shows that
the unique solution $b^*\in(0,1)$ to the above equation is
$b^*=1-\frac1c+\frac1c\log(1-e^{-c}+e^{1-c})$, which can be rewritten as
$$
b^*=b^*(c)=\frac1c\log(1+\frac{e^c-1}e).
$$
Thus, $H(b)$ attains its maximum over $b\in(0,1)$ uniquely at $b^*$.
This proves that $\mathcal{S}(n,M_n^*)$ with $M_n^*$ as in \eqref{optimali} is the asymptotically optimal strategy.
Furthermore, from \eqref{limprobi} and \eqref{H}, the corresponding limiting probability of success
is $\frac{ce^{-c}}{1-e^{-c}}H(b^*(c))$.

We complete the proof of part (i) by  showing  that $\frac{ce^{-c}}{1-e^{-c}}H(b^*(c))=\frac1e$.
From the definition of $b^*(c)$, we have
$$
e^{b^*(c)c}=1+\frac{e^c-1}e,\ \ \ 1-e^{-b^*(c)c}=\frac{e^c-1}{e+e^c-1}.
$$
Therefore,
$$
\begin{aligned}
&H(b^*(c))=(e^{b^*(c)c}-1)\big(1-b^*(c)+\frac1c\log\frac{1-e^{-c}}{1-e^{-b^*(c)c}}\big)=\\
&\frac{e^c-1}e\Big(1-\frac1c\log(1+\frac{e^c-1}e)+\frac1c\log\frac{(1-e^{-c})(e+e^c-1)}{e^c-1}\Big)=\frac{e^c-1}{ec}.
\end{aligned}
$$
Thus, the limiting probability of success
is $\frac{ce^{-c}}{1-e^{-c}}H(b^*(c))=\frac{ce^{-c}}{1-e^{-c}}\frac{e^c-1}{ec}=\frac1e$.

\medskip
\noindent\it Proof of part (ii).\rm\
Let $c>0$ and $\alpha\in(0,1)$ be fixed.  Substitute $q=q_n=1-\frac c{n^\alpha}$ and $M=M_n$  in \eqref{exactform}.
Assume first that $M_n=o(n)$. Since
$$
\frac1{1-q_n^{j-1}}\le \frac{n^\alpha}c,\ \text{for}\ j\ge2,
$$
it follows from \eqref{exactform} that for any $K_1>1$ and any $K_2\in(0,1)$, one has  for sufficiently large $n$,
$$
P_n^{q_n}(\mathcal{S}(n,M_n))\le \frac c{n^\alpha}K_1e^{-cK_2n^{1-\alpha}}\frac{n^\alpha}cn.
$$
Thus $\lim_{n\to\infty}P_n^{q_n}(\mathcal{S}(n,M_n))=0$.

Therefore, from now on we assume that $n=O(M_n)$.
Then it follows from \eqref{exactform} that for any $K_1,K_2>1$, one has for sufficiently large $n$,
$$
P_n^{q_n}(\mathcal{S}(n,M_n))\le \frac c{n^\alpha}K_1e^{-\frac{c(n-M_n-1)}{n^\alpha}}K_2(n-M_n).
$$
The right hand side above converges to 0 as $n\to\infty$ if $\lim_{n\to\infty}\frac{n-M_n}{n^\alpha}=\infty$.
Thus,  from now on, we assume that  $M_n=n-L_n$, with $L_n=O(n^\alpha)$.
Substituting  $M_n=n-L_n$ gives
\begin{equation}\label{optimalii}
\mathcal{P}_n^{q_n}(\mathcal{S}(n,M_n))\sim \frac c{n^\alpha}(1-\frac c{n^\alpha})^{L_n}L_n.
\end{equation}
Let
$$
G_n(L)=(1-\frac c{n^\alpha})^{L}L.
$$
Differentiating, we find that $G_n(L)$ attains its maximum at
$$
L=L_n^*=-\frac1{\log(1-\frac c{n^\alpha})}\sim\frac{n^\alpha}c.
$$
Thus, the asymptotically optimal strategy is $\mathcal{S}(n,M_n^*)$ with $n-M_n\sim L_n^*\sim\frac{n^\alpha}c$.
Substituting $L_n^*$ in \eqref{optimalii} gives
$$
\lim_{n\to\infty}\mathcal{P}_n^{q_n}(\mathcal{S}(n,M_n))=\lim_{n\to\infty}\frac c{n^\alpha}(1-\frac c{n^\alpha})^{\frac{n^\alpha}c}\frac{n^\alpha}c=\frac1e,
$$
which completes the proof of part (ii).
\medskip

\noindent \it Proof of part (iii).\rm\
Fix $q\in(0,1)$. Substituting $M=M_n$ in \eqref{exactform}, one sees that
$\lim_{n\to\infty}\mathcal{P}_n^{q_n}(\mathcal{S}(n,M_n))=0$ if $\lim_{n\to\infty}(n-M_n)=\infty$. Thus, consider
$M_n=n-L$, for an integer  $L\ge1$, and $n$ sufficiently large. Substituting $M_n=n-L$, we have
$$
\lim_{n\to\infty}\mathcal{P}_n^{q_n}(\mathcal{S}(n,M_n))=(1-q)q^{L-1}L.
$$
It is easy to check that the right hand side above attains its maximum over $L\in\mathbb{N}$ at the $L$ for which
$\frac{L-1}L<q\le\frac L{L+1}$. This completes the proof of part (iii).
\hfill $\square$

\end{document}